\documentclass[11pt]{amsart}

\usepackage{enumerate}
\usepackage{amssymb}
\usepackage{graphics}
\usepackage{graphicx}
\usepackage{epsf}

\newtheorem{theorem}{Theorem}[section]
\newtheorem{proposition}[theorem]{Proposition}
\newtheorem{lemma}[theorem]{Lemma}

\newtheorem{definition}[theorem]{Definition}

\theoremstyle{remark}

\newcommand \NN{\mathbf{N}}

\newcommand \eps{\varepsilon}

\newcommand \vs{\vec{s}}
\newcommand \vst{\vs^{\ *}}

\newcommand \bB{\mathbb{B}}
\newcommand \cA{\mathcal{A}}

\newcommand \bC{\mathbf{C}}

\newcommand \restrict{\upharpoonright}

\begin{document}

\title[Cooperative Boolean systems]{Cooperative Boolean systems with generically long attractors II}
\author{Winfried Just
\  and Maciej Malicki
}
\address{Winfried Just: Department of Mathematics, Ohio University, Athens, OH 45701, USA \\
Maciej Malicki: Institute of Mathematics of the Polish Academy of Sciences, Sniadeckich 8, 00-956, Warsaw, Poland}
\email{mathjust@gmail.com, mamalicki@gmail.com}



\begin{abstract}
We prove that cooperativity in Boolean networks precludes a strong notion of sensitive dependence on initial conditions.
Weaker notions of sensitive dependence are shown to be consistent with cooperativity, but if each regulatory functions is binary AND or binary OR, in $N$-dimensional networks they impose an upper bound of $\approx \sqrt{3}^N$ on the lengths of attractors that can be reached from a fraction $p \approx 1$ of initial conditions. The upper bound is shown to be sharp.

\end{abstract}

\maketitle


\section{Introduction}\label{intro}

This paper is a continuation of the investigation in~\cite{JMI} of the extent to which chaos can be generic in cooperative Boolean networks, where \emph{cooperativity} is the absence of any negative interactions whatsoever. The introduction of~\cite{JMI} gives a detailed discussion of our motivation for studying this problem and additional references to related work; here we will repeat only some essentials that are needed for understanding the results of this second part.

Cooperativity is a special case of the more general property of \emph{monotonicity,} which is defined by the absence
of  feedback loops with an odd number of negative interactions.  In monotone flows trajectories converge generically towards an equilibrium under mild regularity hypotheses;  see \emph{e.g.,}  \cite{Enciso:Hirsch:Smith:2008,Hirsch:1988,HalSmith}. In particular, this result implies that chaotic trajectories are not generic in cooperative continuous flows.
Since many natural systems can be modeled with several types of dynamical systems, it is of interest to study whether the above result on nongenericity of chaos in monotone flows has counterparts for Boolean networks.

Chaotic dynamics of Boolean networks are characterized by very long attractors, very few eventually frozen nodes, and high sensitivity to perturbations of initial conditions~\cite{origins}.  These three hallmarks usually, but not always, go together.
Our focus in part~I~\cite{JMI} was on very long attractors.   Since the state space of an $N$-dimensional Boolean network has size~$2^N$, we were interested in upper bounds of the form~$c^N$ for constants $c < 2$.
 As in~\cite{JE, JENDST} we  call an $N$-dimensional Boolean network~\emph{$c$-chaotic} if it does have an attractor of length~$> c^N$. If attractors of this length are reached with probability~$>p$ from a randomly chosen initial condition, then we say that the network is~$p$-$c$-chaotic. Thus $p$-$c$-chaos is a notion of genericity of chaos in terms of very long attractors, and it also implies genericity of chaos in terms of very few eventually frozen nodes (Proposition~3.2. of~\cite{JMI}).

Expected dynamics of so-called \emph{random Boolean networks (RBNs)} tends to become more chaotic as the number of inputs per node increases (see, \emph{e.g.,} the surveys \cite{ACK, Drossel:review, origins}).
The most stringent limitation of this kind is the assumption that the Boolean network is \emph{bi-quadratic,} that is, such that both the number of in- and outputs per node is bounded from above by~2.  The main result of~\cite{JMI} (reproduced below as Theorem~\ref{thm:Ith}) is that that cooperative bi-quadratic Boolean networks can still be $p$-$c$-chaotic.
However, if we require that the system is \emph{strictly bi-quadratic,} that is, all nodes have \emph{exactly} two inputs and two outputs, then even $c$-chaos is possible only for $c <  10^{1/4}$ and the bound is sharp~\cite{JE, JENDST}.
Here we will show that the same bound holds for strictly bi-quadratic $p$-$c$-chaotic cooperative networks (Theorem~\ref{thm:tenroot}).

The main focus of this second part will be the question whether cooperativity limits to some extent the sensitivity to perturbations of initial conditions in Boolean networks.  All by itself, $p$-$c$-chaos does not imply high sensitivity to initial conditions.  In particular, $p$-$c$-chaos can coexist in cooperative Boolean networks, for every $0< p < 1 < c < 2$, with \emph{$p$-coalescence,}  which is the property that  for two randomly chosen initial conditions $\vs(0), \vst(0)$ that differ by  a single-bit flip (have Hamming distance~1) with probability $> p$ there will be some $t>0$ with $\vs(t) = \vst(t)$.

There are several plausible ways of formalizing the notion of sensitive dependence on initial conditions, and we will study three such notions: $p$-instability, which in cooperative Boolean networks is equivalent to the negation of $p$-coalescence, $p$-$D$-decoherence, and $p$-$\alpha$-$q$-decoherence.  Intuitively, the latter two notions mean that a single-bit perturbation to a randomly chosen initial condition will with high probability lead to trajectories that have a relatively large Hamming distance infinitely often.  It turns out that the strongest of these notions, $p$-$\alpha$-$q$-decoherence, does occur in some $p$-$c$-chaotic Boolean networks (Proposition~\ref{prop:unccopdecohere}), but is precluded by cooperativity (Theorem~\ref{thm:toostrongdamth}).  Thus an analogue of the above mentioned theorem for monotone flows does hold for this particular formalization of the notion of chaos in Boolean networks.

For $p$-instability and $p$-$D$-decoherence the situation is more subtle. We will show that for every $0 < p < 1 < c < 2$ there are bi-quadratic cooperative Boolean networks that are simultaneously $p$-unstable and $p$-$c$-chaotic (Theorem~\ref{thm:cppth}). But strictly bi-quadratic such networks can exist only if $c < \sqrt{3}$
(see Subsection~\ref{subsec:strict-quadratic}), and  we show that this bound is again optimal (Theorem~\ref{thm:sqrt3}).  In general, cooperative $p$-$c$-chaotic Boolean networks can exhibit  arbitrarily strong forms of $p$-$D$-decoherence (Theorem~\ref{thm:2dam-p}).  We prove that some versions of this property can occur under the additional assumptions that the network is bi-quadratic (Theorem~\ref{thm:decoh}) or even strictly bi-quadratic (Theorem~\ref{thm:decoh-bi}).  It remains an open problem to determine the maximal amount of $p$-$D$-decoherence that is possible under these additional assumptions.

\section{Terminology}\label{sec:terminology}

Our terminology will be the same as in~\cite{JMI}, where the reader can find all relevant definitions.  Here we will only clarify some key points that are crucial for understanding the formulation of our results.

The symbol~$[N]$ denotes the set $\{1, \ldots , N\}$, which is also the domain of $N$-dimensional Boolean vectors~$\vs \in 2^{[N]}$.  Each such $\vs = (s_1, \ldots , s_N)$ is the characteristic function of the set $A_{\vs} \subseteq [N] = \{i \in [N] : \, s_i =1\}$. It will sometimes be convenient to work with sets $A_{\vs}$ instead of Boolean vectors~$\vs$.  Note that in this interpretation a Boolean function~$f$ is cooperative, which can be defined as preserving the coordinatewise partial order, if and only if it preserves the  subset relation, that is, $A_{\vs} \subseteq A_{\vst}$ implies $A_{f(\vs)} \subseteq A_{f(\vst)}$.  It follows that every partial Boolean function on a set of pairwise incomparable Boolean vectors can be extended to a cooperative total Boolean function (see Proposition 2.1 of~\cite{JMI}); a fact that we will use several times.

The symbol $|\vs|$ denotes the number of coordinates~$i$ with~$s_i = 1$; equivalently, $|\vs| = |A_{\vs}|$.

The Hamming distance $H(\vs, \vst)$ between two Boolean vectors \\ $\vs = (s_1, \ldots , s_N)$ and $\vst = (s_1^*, \ldots , s_N^*)$ with the same domain is the number of~$i$ with $s_i \neq s_i^*$. Vectors with a Hamming of~1 are said to differ by a single bit flip.

As in~\cite{JMI} and elsewhere in the literature, we will use the terms `Boolean system' and `Boolean network' interchangeably.  But we will carefully distinguish these dynamical systems from `Boolean circuits' and `Boolean input-output systems' which are layered arrangements of Boolean gates that calculate certain Boolean functions.
Boolean input-output systems, as opposed to Boolean circuits, allow feedback loops between the variables; both structures can be incorporated as building blocks into Boolean networks to achieve desired dynamics.

A Boolean input-output system is bi-quadratic if every of its variables  has indegree and outdegree at most~2, where the indegree of a variable is the number of variables  its regulatory function takes input from and the outdegree is to number of variables for which it serves as input.  The system is strictly bi-quadratic if, in addition, every of its variables except the input variables has indegree exactly~2.  Note that  for Boolean networks the latter requirement already implies that each variable must have outdegree exactly~2 as well (since the sum of indegrees in any directed graph is equal to the sum of outdegrees), but due to external inputs and outputs this implication is in general false for Boolean input-output systems.

\section{Statement of the results}\label{sec:results}

For easier reference, we state the main result of part~I~\cite{JMI}.

\begin{theorem}\label{thm:Ith}
Given any $0 < p < 1$ and $1 < c < 2$, for all sufficiently large $N$ there exist  $p$-$c$-chaotic, $p$-coalescent, $N$-dimensional bi-quadratic cooperative Boolean networks.
\end{theorem}

\subsection{Sensitivity to initial conditions}\label{subsec:sensitivity}

Our first formal definition of high sensitivity to initial conditions is the notion of $p$-instability that was introduced in~\cite{JE}. A Boolean system is \emph{$p$-unstable} if a random single-bit flip in a randomly chosen initial state moves the trajectory into the basin of attraction of a different attractor with probability at least~$p$. Note that  for cooperative Boolean networks $p$-instability is the same as the negation of ($1-p$)-coalescence: If $\vs(0), \vst(0)$ are two initial conditions that differ at exactly one variable, then we must have either $\vs(0) < \vst(0)$ or $\vs(0)> \vst(0)$; wlog assume the former.  Then cooperativity implies that $\vs(t) \leq \vst(t)$ for all times~$t$.  If the inequality is strict for all~$t$, then the two trajectories must reach different attractors, since in cooperative Boolean networks every two states in a given attractor are incomparable (see, \emph{e.g.,} \cite{JE1}). If equality holds for some~$t$, then the two trajectories coalesce.

We will prove

\begin{theorem}\label{thm:cppth}
Given any $0 < p < 1$ and $1 < c < 2$, for all sufficiently large $N$ there exist  $p$-$c$-chaotic and $p$-unstable $N$-dimensional bi-quadratic cooperative Boolean networks.
\end{theorem}

Another hallmark of chaotic dynamics in Boolean networks is extensive damage propagation, which means that a small perturbation (such as a single-bit flip in an initial condition) tends to spread to a significant proportion of the nodes.  The definition of $p$-instability does not account for this phenomenon.  There are a number of possible ways to formally define extensive damage propagation; we will study here two such notions that require a significant proportion of nodes to be affected when the trajectories already have reached their attractors.

\begin{definition}\label{damdef}
Let $D(N)$ be a  function on the set of positive integer.  An $N$-dimensional Boolean network exhibits \emph{$p$-$D$-decoherence} if with probability $\geq p$ a random one-bit flip $\vst(0)$ in a randomly chosen initial condition $\vs(0)$ results in trajectories with the property that $H(\vs(t), \vst(t)) \geq D(N)$ for infinitely many times $t > 0$. In particular, if $D(N) = \alpha N$ for some constant $\alpha > 0$, then we will refer to $p$-$D$-decoherence as   \emph{$p$-$\alpha$-decoherence.}
\end{definition}

Note that $p$-$\alpha$-decoherence means that for infinitely many~$t$ the Hamming distance will be at least a fraction of~$\alpha$ of the size of the state space.  Our next definition requires this to happen also sufficiently frequently.

\begin{definition}\label{damdefq}
A Boolean network exhibits \emph{$p$-$\alpha$-$q$-decoherence} if with probability $\geq p$ a random one-bit flip $\vst(0)$ in a randomly chosen initial condition $\vs(0)$ results in trajectories with the property that for all sufficiently large $t^* > 0$, the proportion of times
$t \in [0, t^*]$ for which the Hamming distance satisfies
$H(\vs(t), \vst(t)) \geq \alpha N$ is at least~$q$.
\end{definition}

Note that $p$-$\alpha$-$q$-decoherence implies $p$-$\alpha$-decoherence, which in turn implies the negation of ($1-p$)-coalescence, that is, $p$-instability.
Thus in a sense, $p$-$\alpha$-$q$-decoherence is the strongest possible form of  sensitivity to initial conditions.  It turns out that this notion is still consistent
with $p$-$c$-chaos in general, but not with cooperativity.

\begin{proposition}\label{prop:unccopdecohere}
Let $0 < \alpha, p, q < 1$ and $1 < c < 2$.  For all  sufficiently large~$N$ there exist $N$-dimensional $p$-$c$-chaotic Boolean networks that are $p$-$\alpha$-$q$-decoherent.
\end{proposition}

\begin{theorem}\label{thm:toostrongdamth}
For every $\alpha > 0$ and $0 < p < 1$  there exists $N_{\alpha, p}$ such that no cooperative Boolean network of dimension~$N \geq N_{\alpha, p}$ can have the property that for some fixed time $t > 0$ with  probability $\geq p$ a single-bit flip in a randomly chosen initial condition leads to trajectories with
$H(\vs(t), \vst(t)) \geq \alpha N$.  In particular, for any $q > 0$, no cooperative Boolean network of sufficiently large dimension can
exhibit
$p$-$\alpha$-$q$-decoherence.
\end{theorem}

  Thus
$p$-$\alpha$-$q$-decoherence is a chaos-like property of the dynamics that is precluded by cooperativity.
In contrast, the weaker property of $p$-$\alpha$-decoherence is consistent with cooperativity  and $p$-$c$-chaos at the same time.

\begin{theorem}\label{thm:2dam-p}
Let $0 < \alpha, p < 1 < c < 2$. Then for all sufficiently large $N$  there exist  $N$-dimensional  cooperative Boolean networks that are $p$-$c$-chaotic and exhibit
$p$-$\alpha$-decoherence.
\end{theorem}

The networks constructed in our proof of Theorem~\ref{thm:2dam-p} are not subject to any limitations on the number of inputs or outputs per variable and it is of interest to investigate how much damage progation is possible in bi-quadratic cooperative Boolean networks.  We will give a proof of the following result.

\begin{theorem}\label{thm:decoh}
Let $0 < \alpha < 0.5$ and  $0 < p < 1-2\alpha < 1 < c < 2^{1/(1-2\alpha)}$. Then for all sufficiently large $N$  there exist  $N$-dimensional  cooperative bi-quadratic Boolean networks that are $p$-$c$-chaotic and exhibit
$p$-$\alpha$-decoherence.
\end{theorem}

While we do not know whether the bounds on $\alpha, p$, and~$c$ in Theorem~\ref{thm:decoh} are optimal, we conjecture that there are some nontrivial bounds on these parameters in bi-quadratic cooperative networks, that is, we conjecture that the analogue of Theorem~\ref{thm:2dam-p} fails for this class of Boolean networks.

\subsection{Strictly bi-quadratic networks}\label{subsec:strict-quadratic}

The theorems in \cite{JE} give upper bounds on $c < 2$ for $c$-chaotic, cooperative bi-quadratic Boolean networks that have a fixed positive proportion of strictly quadratic regulatory functions.  In particular, if a network is strictly bi-quadratic, the bound is
$10^{1/4}$ and it can be attained. The question is whether a similar result holds for $p$-instability.  Here  we will prove that the same bound is optimal for $p$-$c$-chaotic Boolean networks, that is, we will prove:

\begin{theorem}\label{thm:tenroot}
Let $0 < p < 1$ and $1 < c < 10^{1/4}$. Then for all sufficiently large $N$ there exist  $p$-$c$-chaotic, $p$-coalescent $N$-dimensional strictly bi-quadratic cooperative Boolean networks.
\end{theorem}

The question  arises how much $p$-$c$-chaos and $p$-instability one can have simultaneously in a strictly bi-quadratic cooperative Boolean network. We will prove the following result:

\begin{theorem}\label{thm:sqrt3}
Let $0 < p < 1 < c < \sqrt{3}$. Then for all sufficiently large $N$ there exist  $p$-$c$-chaotic and $p$-unstable $N$-dimensional strictly bi-quadratic cooperative Boolean networks.
\end{theorem}

Note that $\sqrt{3} < 10^{1/4}$. It turns out that Theorem~\ref{thm:sqrt3} is optimal. In order to formally prove this, let us introduce some new terminology.

Define $q_b(c, p)$ as the supremum of all $q$ such that for all sufficiently large~$N$ there exists a strictly bi-quadratic cooperative $p$-$c$-chaotic $N$-dimensional $q$-unstable Boolean network.

Similarly, define $q(c)$ as the supremum of all $q$ such that for all sufficiently large~$N$ there exists a  cooperative $c$-chaotic $N$-dimensional $q$-unstable Boolean network in which all variables have indegree exactly~2.

\medskip

Since $p$-$c$-chaotic networks are automatically $c$-chaotic, for every $p > 0$ the inequality $q_b(c, p)  \leq q(c)$ holds.
In this terminology  Theorem~\ref{thm:sqrt3} simply says that
$q_b(c, p) = 1$ for all $c < \sqrt{3}$ and~$p < 1$.

On the other hand, Theorem~5 of~\cite{JE}  says that for all $c \leq 2$

\begin{equation}\label{eqn:qcest}
 \sqrt{3} < c < 2 \rightarrow q(c) \leq  0.75 +  \frac{\ln (0.5c)}{2\ln 0.75}.
\end{equation}

Notice that on the interval~$[\sqrt{3}, 2]$ the right-hand side of~(\ref{eqn:qcest}) is a function that strictly decreases from~1 to~0.75.
Since  $q_b(c, p) \leq q(c)$, it follows that Theorem~\ref{thm:sqrt3} is in some sense optimal.

It may be of interest to investigate optimal bounds for $q_b(c, p)$ and related functions if $\sqrt{3} < c < 2$.  We wish to leave this as an open problem.

We also don't know whether $p$-$\alpha$-decoherence is possible at all, for \\ any~$p, \alpha > 0$, in strictly bi-quadratic Boolean networks.  However, a very slight weakening of it is still consistent in such networks, even in the presence of $p$-$c$-chaos.

\begin{theorem}\label{thm:decoh-bi}
Let $0 < p <  1 < c < \sqrt{3}$. Then there exists a constant $\Theta = \Theta(p,c) > 0$ such that for all sufficiently large $N$  there exist  $N$-dimensional  cooperative strictly bi-quadratic Boolean networks that are $p$-$c$-chaotic and exhibit
$p$-$\frac{N}{\Theta\log(N)}$-decoherence.
\end{theorem}

\section{Damage propagation and $p$-instability}\label{sec:damage-propagation}

Here we prove all results that were announced in the previous section that do not require any knowledge of the proof of Theorem~\ref{thm:Ith}. The proofs of  Theorems~\ref{thm:cppth} and~\ref{thm:decoh}--\ref{thm:decoh-bi}  rely to some extent on the construction that was used in~\cite{JMI}  and will be given in the next section.

\bigskip

\noindent
\textbf{Proof of Proposition~\ref{prop:unccopdecohere}:} Fix $\alpha, p, q, c$ as in the assumption. Let $N$ be sufficiently large such that
\begin{equation}\label{eqn:attr-small-enough}
c^N + 1 < \frac{1-p}{N-1} 2^{N-2}.
\end{equation}

It will be convenient for this proof to treat the states of~$\bB$ as subsets of~$[N]$ instead of Boolean vectors. Let $L$ be the  integer that satisfies $c^N < L \leq c^N + 1$.   Choose an indexed set $\cA = \{A_\ell: \ \ell \in [L]\}$ with $A_\ell \subseteq [N-1]$ and
define an updating function~$f$ for~$\bB$ as follows:

\begin{equation}\label{eqn:oscillating-f}
\begin{split}
&f(A_\ell) =  A_{\ell + 1}\ \mbox{for} \ \ell \in [L-1];\\
&f(A_L) =  A_{1};\\
&f([N] \backslash A_\ell) =  [N] \backslash A_{\ell+1}\ \mbox{for} \ \ell \in [L-1];\\
&f([N] \backslash A_L) =  [N] \backslash A_1;\\
&f(B) =  A_{1} \ \mbox{if}\ B \cap [N-1], [N-1] \backslash B  \notin \cA \ \mbox{and} \ |B| \ \mbox{is odd};\\
&f(B) =  [N] \backslash A_1 \ \mbox{if}\ B \cap [N-1], [N-1] \backslash B \notin \cA \ \mbox{and} \ |B| \ \mbox{is even}.
\end{split}
\end{equation}

Note that $A_\ell \neq [N] \backslash A_{\ell'}$ for all $\ell, \ell'$.

Now consider  initial conditions $\vs(0), \vst(0)$ where $\vs(0)$ is randomly chosen and $\vst(0)$ is obtained by a random single-bit flip, and let $B, B^*$ be the sets of indices in~$[N-1]$ with $s_i(0)=1$ and $s^*_i(0) = 1$ respectively.
By~(\ref{eqn:attr-small-enough}), with probability~$>p$, neither of the sets $B \cap [N-1], [N-1] \backslash B, B^* \cap [N-1], [N-1] \backslash B^*$ will be in~$\cA$, and the last two clauses of the definition of the updating function~$f$ apply.   Hence wlog $f(B) = A_{1}$ and $f(B^*) = [N] \backslash A_{1}$.  Thus at time~1 the system will have entered two different attractors of length~$> c^N$ for these initial conditions, and we will have $H(\vs(t), \vst(t)) = N$ for all $t > 0$. $\Box$

\bigskip

\noindent
\textbf{Proof of Theorem~\ref{thm:toostrongdamth}:} Let $\alpha, p$ be as in the assumptions and let $N_{\alpha, p}$ be the smallest positive integer $N$ such that for all $k \in [N]$

\begin{equation}\label{binomest}
\frac{\binom{N}{k}}{2^N} < \frac{p\alpha}{2}.
\end{equation}

Let $\bB$ be a cooperative Boolean system of dimension $N \geq N_{\alpha, p}$. By symmetry we may focus in this argument  on the case where a single bit is flipped from~0 to~1.  Fix $t > 0$ and let $r$ be the probability that a single-bit flip from~0 to~1 in a randomly chosen initial condition leads to trajectories with
$H(\vs(t), \vst(t)) \geq \alpha N$.
Assume towards a contradiction that  $r \geq p$.
For each $k \in \{0, \ldots , N-1\}$, let $p_k$ be the conditional probability that a single-bit flip from $0$ to~$1$ in a randomly chosen initial state~$\vs(0)$
given that~$|\vs(0)|=k$ results in trajectories with

\begin{equation}\label{Hammingeq}
H(\vs(t), \vst(t)) \geq \alpha N.
\end{equation}

Note that in this case $|\vst(0)| = k+1$ and  $\vs(0) < \vst(0)$. Cooperativity implies that

\begin{equation}\label{coopt}
\vs(t) < \vst(t).
\end{equation}

Let $L = \{k: \, p_k \geq \frac{r}{2}\}$ and let $K = \{k: \, p_k < \frac{r}{2}\}$. By~(\ref{binomest}),

\begin{equation}\label{psum}
r \leq \sum_{k = 0}^{N-1} \frac{p_k \binom{N}{k}}{2^N} = \sum_{k \in L} \frac{p_k \binom{N}{k}}{2^N} + \sum_{k \in K} \frac{p_k \binom{N}{k}}{2^N} < \frac{|L|p\alpha}{2} + \frac{r}{2}.
\end{equation}

Under the assumption $r \geq p$ this implies

\begin{equation}\label{manyps}
|L| = |\{k: \, p_k \geq \frac{r}{2}\}| > \frac{1}{\alpha}.
\end{equation}

Now consider a randomly chosen permutation~$\pi$ of~$[N]$, and let $\vs^{\, k,\, \pi}(0)$ be the characteristic function of the set
$\{j : \, \pi(j) < k\}$. Define random variables $X_k$ such that $X_k(\pi)$ takes the value~$1$ if
$H(\vs^{\, k,\, \pi}(t), \vs^{\, k+1,\, \pi}(t)) \geq \alpha N$
and takes the value~$0$ otherwise.  Let $X = \sum_{k=0}^{N-1} X_k$.  Then $E(X_k) = p_k$ for all~$k$ and hence
$E(X) =   \sum_{k=0}^{N-1} p_k$. By~(\ref{manyps}), $E(X) > \frac{1}{\alpha}$, and it follows that there exists at least one permutation~$\pi$
with $X(\pi) > \frac{1}{\alpha}$.  But existence of such a permutation would require in view of~(\ref{Hammingeq}) and~(\ref{coopt}) that there exist initial states $\vs^{\, 0, \, \pi }(0) < \vs^{\, 1,\, \pi}(0) < \dots < \vs^{\, J,\, \pi}(0)$ with $J > \frac{1}{\alpha}$ such that
$\vs^{\, 0, \, \pi}(t) < \vs^{\, 1, \, \pi}(t) < \dots < \vs^{\, J, \, \pi}(t)$ are characteristic functions of sets $A_j$ with $A_j \subset A_{j+1} \subseteq [N]$ and
$|A_{j+1} \backslash A_j| \geq \alpha N$, which leads to a contradiction.

It remains to show how the first part of the theorem implies the second. Fix $\alpha, p, q$ as in the definition of $p$-$\alpha$-$q$-decoherence.
 For each $t \geq 0$ consider the random variable $\xi_t$ on the space of all pairs $(\vs(0), \vst(0))$ that result from a random bit flip in an initial condition that takes the value~$1$ if $H(\vs(t), \vst(t)) \geq \alpha N$ and takes the value~$0$ otherwise. The first part of the proof shows that as long as
 $N \geq N_{\alpha, pq}$ we will have

 \begin{equation}\label{eqn:Exit}
 E(\xi_t) = P(\xi_t = 1) < pq.
 \end{equation}

Now fix $t^* \geq 0$ and let $\eta = \sum_{t=0}^{t^*} \xi_t$.  If~(\ref{eqn:Exit}) holds, then

\begin{equation*}\label{eqn:Eeta}
 q t^* P(\eta \geq qt^*) \leq  E(\eta)  < pq t^*,
 \end{equation*}
and it follows that

$$P(\eta \geq qt^*)    < p,$$
 which contradicts $p$-$\alpha$-$q$-decoherence. $\Box$

\bigskip

\noindent
\textbf{Proof of Theorem~\ref{thm:2dam-p}:} Let $\alpha, p, c$ be as in the assumptions.  Fix the smallest positive integer~$z$ with
$p < 1 - 2^{-z + 2}$, and fix $\gamma > 0$ and $N_\gamma > 2z$ such that the following inequality holds for all $N > N_\gamma$:

\begin{equation}\label{gammaeq}
\sum_{k = \lceil N/2 -  \gamma \sqrt{N}\rceil + 1}^{k = \lceil N/2 +  \gamma \sqrt{N}\rceil - 1} \frac{\binom{N}{k}}{2^N} >
p + 2^{-z + 2}.
\end{equation}

For $N > N_\gamma$, let $w: = \lceil N/2 -  \gamma \sqrt{N}\rceil$ and $u := \lceil N/2 +  \gamma \sqrt{N}\rceil$. We will assume for sake of simplicity that
$u-w$ is even.

By assumption, $[2z] \subset [N]$.  The Boolean variables~$s_i$ with~$i \in [2z]$ will play a special role in controlling cooperativity of the Boolean system that we are going to construct.

Let~$Z$ be the set of all states~$\vs$ that satisfy the following conditions:

\begin{equation}\label{zeroineq}
\begin{split}
&\exists i \in [z] \, \exists j \in [2z] \backslash [z] \ s_i = 1 \ \& \ s_j = 0,\\
&w \leq |\vs| \leq u.
\end{split}
\end{equation}

Let $N > N_\gamma$ and consider a randomly chosen initial condition $\vs(0)$ and any condition~$\vst(0)$ obtained from it by a single-bit flip.
The probability that the first line of~(\ref{zeroineq}) fails for $\vs(0)$ or~$\vst(0)$ is less than~$2^{-z + 2}$, and~(\ref{gammaeq}) implies that $P(w+1 \leq |\vs(0)| \leq u-1) > p + 2^{-z + 2}$. It follows that with probability  $> p$ both  $\vs(0), \vst(0) \in Z$.

We will construct systems~$\bB$ of dimension $N > N_\gamma$ as follows. Let $J = u - w$. For each $j \in J$ we will
specify a periodic orbit $A_j = \{\vs^{\, j}(i): \, i \in [L]\}$ of length $L > c^N$, where  $\vs^{\, j}(i+1)$ is the successor state in~$\bB$ of $\vs^{\, j}(i)$ for all $i < L$,  in such a way that

\smallskip

\noindent
(i) $\vs^{\, j}(i) < \vs^{\, j+1}(i)$ for all $j \in [J]$ and $i \in [L]$,

\smallskip

\noindent
(ii) $\vs^{\, j}(i)\restrict [z] = \vec{0}$ and $\vs^{\, j}(i)\restrict [2z] \backslash [z] = \vec{1}$  for all $j \in [J]$ and $i \in [L]$,

\smallskip

\noindent
(iii) $H(\vs^{\, j}(j), \vs^{\, j+1}(j)) \geq \alpha N$ for all $j \in [J]$, and

\smallskip

\noindent
(iv) for $i \neq i'$ and any $j, j' \in J$ the states  $\vs^{\, j}(i)$ and $\vs^{\, j'}(i')$ are incomparable with respect to the coordinatewise partial order.

\medskip

This part of the construction defines a partial Boolean updating function~$f$ of~$\bB$ on the set
$A = \bigcup_{j \in [J]} A_j$.  By~(i), (ii) and~(iv) this function is cooperative.

 Note that by~(ii), all states in~$A$ are incomparable with all states in~$Z$.  Thus if we define the restriction of the Boolean updating function~$f$ to~$Z$ so that it is cooperative, then automatically $f\restrict (A \cup Z)$ will be cooperative. Now consider~$\vs \in Z$.  Then for a unique~$j \in [J]$ we have $|\vs| = w + j$, and we define~$f(\vs) = \vs^{\, j}(1)$.  By~(i), this construction results in a cooperative Boolean function on~$Z$. Having defined a cooperative partial Boolean function~$f\restrict (A \cup Z)$  we can extend it by
Proposition~2.1 of~\cite{JMI} to a
cooperative updating function~$f$ on the whole state space~$2^N$ of~$\bB$.

Now consider a random initial condition~$\vs(0)$ and let~$\vst(0)$ be obtained by some one-bit flip in~$\vs(0)$.  Then with probability $> p$
both
$\vs(0), \vst(0) \in Z$, and it follows that there are $j, j' \in [J]$ with $|j-j'| =1$ such that $\vs(1) = \vs^{\, j}(1)$ and
$\vst(1) = \vs^{\, j'}(1)$. Wlog $j' = j+1$  and condition~(iii) implies that $H(\vs(t), \vst(t)) \geq \alpha N$ for infinitely many~$t$, which gives $p$-$\alpha$-decoherence.

Note that when $\vs(0) \in Z$, the trajectory of~$\vs(0)$ will reach one of the attractors $A_j$.  In particular,~(iv) implies that all $s^j(i)$ are pairwise distinct for different~$i$ and fixed~$j$, thus $A_j$ has length~$L > c^n$, and we  get $p$-$c$-chaos  as well.

\medskip

It remains to prove that for sufficiently large~$N$ we can construct a family $\cA = \{A_j: \, j \in [J]\}$ that satisfies conditions~(i)--(iv).

We need $N_0 > N_\gamma$ sufficiently large so that for $N > N_0$ we have

\begin{equation}\label{eqn:N0-first-condition}
c^N < \binom{N - 2(z + u - w)}{\lceil N/2\rceil - (z + u - w)},
\end{equation}

\begin{equation}\label{eqn:N0-second-condition}
\binom{u-w}{(u - w)/2} > u - w, \ \mbox{and}
\end{equation}

\begin{equation}\label{eqn:N0-third-condition}
\frac{N - 2(z + u - w)}{N} > \alpha.
\end{equation}

Conditions~(\ref{eqn:N0-second-condition}) and~(\ref{eqn:N0-third-condition}) will be quite obviously satisfied for all sufficiently large~$N$; condition~(\ref{eqn:N0-first-condition}) follows from the fact that we can make $2^{N - 2(z + u - w)}$ larger than $d^N$ for any $d < 2$ and that
$\binom{2K}{K} \sim \frac{2^{2K}}{\sqrt{K}}$ as $K \rightarrow \infty$.

Fix $N > N_0$.
Let $U, W$ be disjoint subsets of $[N] \backslash [2z]$ such that $|U| = |W| = u - w$, and let $\{a_j: \, j \in \{0\} \cup [u-w]\}$ be a family of pairwise
incomparable subsets of $U$.  We can form this family from subsets of~$U$ of size $(u-v)/2$ each; condition~(\ref{eqn:N0-second-condition}) implies that there will be enough such sets to choose from. Similarly, by~(\ref{eqn:N0-first-condition}) we can choose a family $C = \{c_i: \, i \in [L]\}$ of size~$L > c^N$ of
subsets of $[N] \backslash ([2z] \cup U \cup W)$ that have size~$\lceil N/2\rceil - (z + u - w)$ each and thus are pairwise incomparable.   Let $W_j$ for $j \in [u - w]$ be subsets of~$W$ such that $W_j$ is a proper subset of $W_{j+1}$ for all relevant~$j$.

Then define $\vs^{\, j}(i)$ as the characteristic function of the set

\begin{itemize}
\item $a_i \cup W_j \cup  ([N] \backslash ([z] \cup U \cup W))$ if $i < j  \leq
u - w$,
\item $a_i \cup W_j \cup [2z] \backslash [z]$ if $j \leq i  \leq
u - w$, and
\item $a_0 \cup W_j \cup c_i \cup [2z] \backslash [z]$ if $i >
u - w$.
\end{itemize}

It is straightforward to verify that  conditions~(i)--(iv) hold, with the all-important condition~(iii) following from~(\ref{eqn:N0-third-condition}). $\Box$

\section{Proofs of Theorem~\ref{thm:cppth} and \ref{thm:decoh}--\ref{thm:decoh-bi}}\label{sec:constructs}

The proofs of these theorems are based on the construction that was used in~\cite{JMI} for the proof of Theorem~\ref{thm:Ith} and we will need to review it here to some extent.

Let $0 < p < 1 < c < 2$ be as in the assumptions of Theorem~\ref{thm:Ith}.
In the proof we constructed, for sufficiently large~$N$, a suitable
updating function~$f$ for  Boolean systems $\bB = (2^N, f)$ such that~$f$ was cooperative, bi-quadratic, and worked as required.
The set of Boolean variables $[N]$ was partitioned into a disjoint union $[N] = X \cup Y$, where the set~$X$ in turn was a union of  pairwise disjoint sets~$X_i$, indexed by $i \in I = \{0, 1, \ldots , |I|-1\}$, and all of the same size~$m \leq |I|$. Both $m$ and $|I|$ scale like $\sqrt{N}$.
We singled out some $i_2 > i_1 > i_0 \in I$ and conceptualized the collection of all sets~$X_i$ as a circular data tape, with $f$ simply copying the vector $s_{X_{i+1}}(t)$ to $s_{X_{i}}(t+1)$ for all indices~$i$ with the exception of~$i \in \{i_0, i_1\}$, and also copying $s_{X_{0}}(t)$ to $s_{X_{|I|-1}}(t+1)$. The vectors $s_{X_{i_1}}(t+1), s_{X_{i_0}}(t+1)$  were outputs of special Boolean circuits $B_2, B_1$, which also  took a second input from another Boolean input-output system $B_3$.

A schematic view of the construction is given in Figure~1.

\begin{figure}
 \centering
 \includegraphics[width=5in,height=4in]{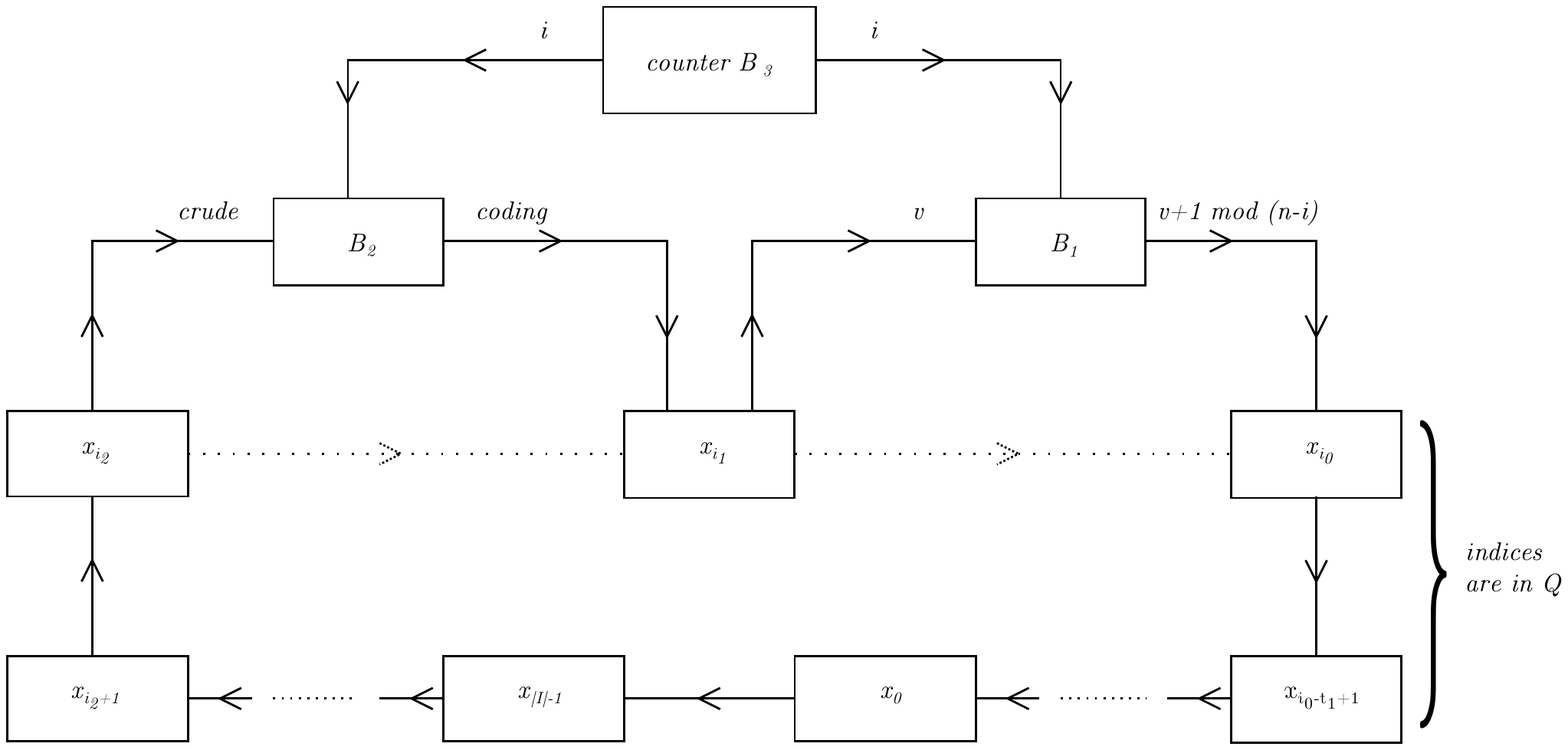}
 \caption{\label{fig:overall}  A schematic view of the construction.}
\end{figure}

Although the figure suggests $0 \notin Q$, it will be convenient here to assume that $i_0 - t_1 + 1 = 0$, so that $Q = \{0, 1, \ldots i_0\}$. This amounts to a circular shift in the indexing and does not alter the construction in any way.

Our proofs will rely on the following facts:

\begin{itemize}
\item[(P1)] For sufficiently large~$N$ it is possible to construct these objects so that the resulting Boolean system will be cooperative, bi-quadratic, and exhibit $p$-$c$-chaos.
\item[(P2)] For all $i \in I \backslash \{i_0, i_1\}$ and all $t$ we will have $|s_{X_{i+1}}(t)| = |s_{X_{i}}(t+1)|$.
\item[(P3)] Each variable in $X_{i_2 + 1}$ has exactly one output variable in the system.
\item[(P4)] With probability $>p$ the trajectory of a randomly chosen initial condition will have the property that for every  time $t = k|I|$ such that $k$ is a \emph{positive} integer, we have $|s_{X_i}(t)| = |X_i|/2$ for all $i \in \{i_2+1, i_2+2, \ldots , |I|-1\}$. This follows from the choice of coding vectors in~\cite{JMI} and the workings of~$B_1, B_2, B_3$.
\item[(P5)] The ratio $|Y \cup \bigcup_{i \leq i_2} X_i|/N$ approaches zero as $N \rightarrow \infty$.
\end{itemize}

 The key requirement that makes this construction work for obtaining $p$-$c$-chaos is the following:

\begin{itemize}
\item[(P6)] With probability $> p$ the following will hold for a randomly chosen initial condition:
if $t$ is any time of the form~$k|I|$, where $k$ is a positive integer, then
with the possible exception of indices~$i \in Q$, each vector $s_{X_i}(t)$ will be a \emph{coding vector,} that is, will code an integer
$v_i(t) \in \{0, \ldots , n-1\}$ for some suitable value of~$n$ that depends on~$N$.  Moreover, again with the possible exception of $i \in Q$, the function~$f$ computes addition of~1 modulo $n-i$  on input $X_i$ and writes the output to $X_i$ after~$|I|$ steps.  Formally, the latter means that for times~$t$ as above

\begin{equation}\label{modaddI0}
\forall \, i \in I \backslash Q \quad v_i(t+ |I|) = v_i(t) + 1 \ mod \ (n - i).
\end{equation}
\end{itemize}

In order to make (P6) work, we will need a suitable~$n$ and, for each~$X_i$, a set of coding vectors~$\bC_i \subseteq 2^{X_i}$.  The choice of the sets $\bC_i$ will be restricted (in the sense of~(P7a) below) by the particular coding scheme that we adopt in a given proof. For technical reasons the construction in~\cite{JMI} needs

\begin{equation}\label{eqn:size-of-m}
m = |X_i| = (1+\eps)\log n \quad \mbox{and} \quad |I| = \beta\log(n),
\end{equation}
where $\eps >0$  and $\beta$ is a positive integer that does not depend on~$n$.  Note that~(\ref{eqn:size-of-m}) implies the scaling laws $m = O(\sqrt{N})$ and  $|I| = O(\sqrt{N})$.   Only the following additional properties of the coding are needed to obtain $p$-$c$-chaotic systems:

\begin{itemize}
\item[(P7)] There are positive integers~$k, \ell$ with $k$ depending only on~$c$ such that
\begin{itemize}
\item[(P7a)]  $\bC_i \subset (C_k)^\ell$, where  $C_k$ is a set of Boolean vectors from $2^k$,  called the \emph{coding scheme,} such that exactly half of their coordinates are $1$'s (so the other half are $0$'s).
\item[(P7b)] The pair $(k, \eps)$ is \emph{$c$-friendly under the given coding scheme,} which means that $\eps$ is rational, $\frac{k}{1+\eps}$ is an integer,

\begin{equation}\label{eqn:logceps1}
{\log(c)} (1 + \eps) <  1, \ \mbox{ and}
\end{equation}

\begin{equation}\label{eqn:friend1eq}
|C_k| \geq 2^{k/(1+\eps)}.
\end{equation}
\end{itemize}
\end{itemize}

A few remarks are in order here.  Of course, the sets of variables $X_i$ are pairwise disjoint, so we cannot literally
make each $\bC_i$ a subset of $(C_k)^\ell$; formally we will need disjoint copies of  $(C_k)^\ell$.  However, we suppress the additional parameter to reduce clutter in our notation. The outputs of the Boolean input-output system~$B_3$ will also code for integers and satisfy property~(P7), so our modifications of the coding scheme will apply to them as well.  The wording chosen in~\cite{JMI} suggests that $C_k$ is the set of \emph{all} Boolean vectors from $2^k$ such that exactly half of their coordinates are $1$s, but this property was never actually used in the proof. Only~(\ref{eqn:size-of-m}) and property~(P7) (which  are taken from Section~8 of~\cite{JMI})
are ever referred to in any part of the construction. In fact, since we need exactly~$n$ codes for nonnegative integers, for most~$N$ not all vectors in $(C_k)^\ell$ are used even in~\cite{JMI} as actual codes.  This makes it possible for us to alter the definition of $C_k$  to more restrictive coding schemes that  will work for our purposes here.

For the description below, it will be convenient to consider a partition of each $X_i$ into pairwise disjoint subsets~$x^j_i$ of size~$k$ each that correspond to the domains of the vectors in~$C_k$ when~$s_{X_i}$ is coding.

The system works as follows: With probability arbitrarily close to~1, for each~$i$ the vector $s_{X_i}(0)$ will be \emph{crude}, which means that there will be $j,{j'}$ such that the restriction of $s_{X_i}(0)$ to $x_i^j$ will take the constant value~0  and the restriction of $s_{X_i}(0)$ to $x_i^{j'}$ will take the constant value~1 (Lemma~9.1 of~\cite{JMI}).  For $i \notin \{i_0, i_0 + 1, \ldots , i_2-1\}$ these crude vectors will be eventually copied to $X_{i_2}$, where they become inputs of the Boolean circuit~$B_2$, which eventually outputs a code for~$0$ to~$X_{i_1}$ for each crude input.  When given an input that is a coding vector, $B_2$ outputs an identical copy of its input.  The vector~$s_{X_{i_2}}$ becomes one of the inputs of the Boolean circuit~$B_1$, which eventually outputs a code for $v + 1  \ mod \ (n - w)$ to $X_{i_0}$ if its input from~$X_{i_2}$ codes the integer~$v$ and the other input that it receives from~$B_3$ codes an integer~$w$. With probability arbitrarily close to~1, the Boolean input-output system~$B_3$ will deliver the correct sequence of inputs to~$B_1$ so that~(\ref{modaddI0}) of property~(P6) holds.  The set~$Q$ indexes those $X_i$ for which the output of~$B_1$ may not be coding, due to the time lag in the calculations of $B_1, B_2, B_3$.

\bigskip

\noindent
\textbf{Proof of Theorem~\ref{thm:tenroot}:}  We need to turn the networks constructed  in~\cite{JMI} into strictly bi-quadratic ones.  The main problem is that in the original construction each vector~$s_{X_i}(t+1)$ was supposed to be a copy of $s_{X_{i+1}}(t+1)$ as long as $i \notin \{i_0, i_0 + 1, \ldots , i_2-1\}$.  This can be easily accomplished by a Boolean circuit~$B_c$ of depth~1 with input variables $X_{i+1}$ and $X_{i}$ as output variables that uses only COPY functions, but there is no analogous strictly bi-quadratic cooperative Boolean circuit.  Fortunately, as can be seen from the above description,  it is not actually necessary that $B_c$ outputs identical copies of \emph{all} possible inputs; it suffices that it does so whenever its input~$s_{X_{i+1}}$ is a coding vector.  Since we may wlog assume that $|I|$ is even, it even suffices to require that $s_{X_i}(t+2) = s_{X_{i+2}}(t)$ holds for all relevant indices~$i$ whenever~$s_{X_{i+2}}(t)$ is coding.  It turns out that there are strictly bi-quadratic cooperative Boolean circuits that work in this way for sets of coding vectors that satisfy property~(P7) as long as $c < 10^{1/4}$.

We will also need that~$B_c$ outputs crude vectors for crude inputs, so that $B_2$ will receive the kind of inputs that are expected for random initial conditions.
However, it follows immediately from the definition of crudeness that this will be automatically satisfied if $B_c$ is cooperative, strictly quadratic (thus uses only
binary AND and OR gates), and such that all inputs for variables in $x^j_i$ reside in~$x^j_{i+1}$, as will be the case in our constructions.

Let us now present two coding schemes that allow for implementation of this idea.  The first one will only be used in the proofs of some of our other theorems, but it is easier to understand and we want to describe it here as an illustration.   Let $X_i, X_{i+1}$ be consecutively enumerated by $\kappa(i, \lambda), \kappa(i+1, \lambda)$, where
$\lambda \in [m]$.  Let $C_k$ be the set of Boolean vectors $\vs \in 2^k$ such that $s_{2\kappa - 1} \leq s_{2\kappa}$ for all $\kappa \in [k/2]$ and $|\vs| = k/2$, and let $\bC_i$ be the corresponding sets of codes in the sense of~(P7a). Define $B_{c,r} = B_{c,r}(X_{i+1}, X_i)$ as the Boolean circuit of depth~1 with input~$X_{i+1}$, output $X_{i}$ and regulatory functions defined by

\begin{equation}\label{eqn:regulatory-Bcr}
\begin{split}
s_{\kappa(i, 2\lambda-1)}(t+1) &= s_{\kappa(i+1, 2\lambda-1)}(t) \wedge s_{\kappa(i+1, 2\lambda)}(t) \quad \mbox{for} \ \lambda \in [m/2];\\
s_{\kappa(i, 2\lambda)}(t+1) &= s_{\kappa(i+1, 2\lambda-1)}(t) \vee s_{\kappa(i+1, 2\lambda)}(t) \quad \mbox{for} \ \lambda \in [m/2].
\end{split}
\end{equation}

Now it is clear from~(\ref{eqn:regulatory-Bcr}) that the Boolean circuit $B_{c,r}$ is strictly bi-quadratic and, as long as the sets $x^j_i$ and~$x^j_{i+1}$ are consecutive intervals in $X_{i}, X_{i+1}$,  it will copy $\vs \in \bC_{i+1}$ to its counterpart in~$\bC_i$.  Moreover, it will map crude vectors to crude ones.
The circuit~$B_{c,r}$ has the additional useful property of mapping each $\vs$ to $f(\vs)$ such that $|\vs| = |f(\vs)|$, regardless of whether or not~$\vs$ is coding. This will allow us to retain property~(P2).
For this reason we will call the coding scheme that we just defined \emph{robust coding.}

\begin{lemma}\label{lem:friendlem1}
Suppose $1 < c < \sqrt{3}$.  Then there exist a rational $\eps = \eps(c) > 0$ and a positive even integer $k = k(c)$ such that
 the pair $(k, \eps)$ is $c$-friendly under robust coding.
\end{lemma}

\noindent
\textbf{Proof:} Let $\eps$ be rational, such that $\frac{k}{1+\eps}$ is an integer,
 and

\begin{equation}\label{eqn:eps-for-robust}
\log(c) < \frac{1}{1+\eps} < \log(\sqrt{3}).
\end{equation}

Such $\eps$ exists by our assumption on~$c$, and~(\ref{eqn:eps-for-robust}) implies~(\ref{eqn:logceps1}).

Fix an even integer~$k$.  We can think about the vectors $\vs \in C_k$ as outcomes of the experiment of randomly and independently drawing
$(s_{2\kappa - 1} \leq s_{2\kappa})$ from the set $\{(00), (01), (11)\}$ with the uniform distribution.  Then $|\vs|$ is a random variable with mean $E(\vs) = k/2$, and the space of all possible outcomes has size $3^{k/2}$.  The Central Limit Theorem implies that the probability of obtaining a vector in $C_k$, that is, an outcome with $|\vs| = k/2$,  scales like
$k^{-0.5}$. Thus  by~(\ref{eqn:eps-for-robust})  we have for some constant~$\rho > 0$ and for sufficiently large $k$

\begin{equation}\label{eqn:size-Ck-robust}
|C_k| \geq \rho 3^{k/2 - 0.5\log(k)}  > 2^{k/(1+\eps)},
\end{equation}
and~(\ref{eqn:friend1eq}) follows. $\Box$

\bigskip

The coding scheme that we will actually use in the proof of Theorem~\ref{thm:tenroot} is more complicated and we will refer to it as \emph{subtle coding.} The relevant sets $C_k$ and the corresponding Boolean
circuit~$B_{c,s}(X_{i+1}, X_i)$ were implicitly described in detail in Section 4.2 of~\cite{JE} and Section~3 of~\cite{JENDST}, and we refer the reader to these publications for details.  Here we only want to record  the key properties that will be used in our arguments.

\begin{lemma}\label{lem:friendlem2}
(a) Suppose $1 < c < 10^{1/4}$.  Then there exist a rational $\eps = \eps(c) > 0$ and a positive even integer $k = k(c)$ such that
 the pair $(k, \eps)$ is $c$-friendly  under subtle coding.

\smallskip

\noindent
(b) There exists a strictly bi-quadratic Boolean circuit $B_{c,s}(X_{i+1}, X_i)$ of depth~1 such that if $s_{X_i}$ is a coding vector under the subtle coding scheme, then the concatenation $B_{c,s}(X_{i+1}, X_i) \circ B_{c,s}(X_{i+2}, X_{i+1}) $ outputs an identical copy of $s_{X_{i+2}}$ after two steps and writes it to the variables
in~$X_i$.
\end{lemma}

\noindent
\textbf{Proof:} For the proof of part~(a), choose~$c_1$ with $c < c_1 < 10^{1/4}$. Lemma~3.1 of~\cite{JENDST} states (in a slightly different terminology) that for sufficiently large~$k$ that are divisible by~8 we will have $|C_k| > c_1^k$. Now let $\eps$ be a rational such that $\frac{k}{1+\eps}$ is an integer,
 and

\begin{equation}\label{eqn:eps-for-subtle}
\log(c) < \frac{1}{1+\eps} < log(c_1).
\end{equation}

Then

\begin{equation}\label{eqn:size-Ck-subtle}
|C_k| \geq \rho c_1^k > 2^{k/(1+\eps)},
\end{equation}
as required.

For the proof of part~(b) we refer the reader to~\cite{JE} or~\cite{JENDST}. $\Box$

\bigskip

Unfortunately, subtle coding does not preserve $|\vs|$.  Thus in the proof of Theorem~\ref{thm:sqrt3} we will use robust coding instead.  Moreover, the output vector $f(\vs)$ of $B_{c,s}$ is not usually an exact copy of~$\vs$ even if $\vs$ is coding.  However, applying the operation twice produces a copy~$f\circ f(\vs)$ of $\vs$ whenever $\vs$ is a subtle code.  This is sufficient for our purposes; as we already mentioned above, we only need that $s_{X_i}(t+2)$ is a copy of $s_{X_{i+2}}(t)$ for all relevant $i$ and~$t$.

Lemmas~\ref{lem:friendlem1} and~\ref{lem:friendlem2}  imply that instead of straight copying in the construction of~\cite{JMI} we can use the circuits $B_{c,r}$ with robust coding for all $0 < p < 1 < c < \sqrt{3}$ or $B_{c,s}$ with subtle coding for all $0 < p < 1 < c < 10^{1/4}$. This will not affect the other technical arguments  of the construction in~\cite{JMI} and give us
$p$-$c$-chaotic systems. Alas, it will not give us all by itself strictly bi-quadratic networks. We need to alter the construction in such a way that \emph{every} variable has exactly two inputs and exactly two outputs.  If we use robust or subtle coding, then this will be true for the variables in the relevant sets~$X_i$, but not automatically for the variables in the remaining parts of the system, in particular, for the variables in $B_1, B_2, B_3$.  We need a tool for adding redundant inputs to some variables that will not substantially alter the long-range dynamics of the whole system. The following lemma gives us such a tool.

\begin{lemma}\label{lem:ones}
For every $0 < q < 1$  there exists a Boolean system $B_q$ of depth $d= \lceil \log(- \log(1-q))\rceil$ with $< 2\lceil - \log(1-q)\rceil$ variables that satisfies the following. Except for one variable $i_q$ whose indegree is $1$, the indegree of every variable is $2$, and, except for one variable $o_q$ whose outdegree is $1$, the outdegree of every variable is $2$. Moreover, with probability $\geq q$ the value of $o_q$ will be~1 at all times $t \geq d$, regardless of the initial values of all the variables in $B_q$, and the trajectory of~$i_q$.
\end{lemma}

\noindent
\textbf{Proof:} Let $q$ be as in the assumption and let $d$ be as in the statement of the Lemma.  It will be convenient to let the variables of the system be binary sequences~$\sigma$ of length at most~$d$, where $o_q$ is the empty sequence, and $i_q$ is the zero sequence of length $d$. The sequences of length $d$ will constitute the lowest level~1 of the  variables of the system.  For $\sigma$ of length $< d$ we let

\begin{equation}\label{eqn:ones-reg-higher}
s_{\sigma}(t+1) = s_{\sigma^\frown0}(t) \vee s_{\sigma^\frown1}(t).
\end{equation}

It is easy to see that the total number of variables on levels  $> 1$ is $2^d-1$. Therefore, there exists a bijection $\varphi$ between all the variables on level  $1$ except for $i_q$, and the variables on higher levels. For variables $\sigma$ on level~$1$ that are distinct from $i_q$ we define

\begin{equation}\label{eqn:ones-reg-lowest}
s_{\sigma}(t+1) = s_{\sigma}(t) \vee s_{\varphi(\sigma)}(t),
\end{equation}
and we let $s_{i_q}(t+1)=s_{i_q}(t)$.

Thus the variable $o_q$ will take the value~$0$ at time~$d$ only if $s_{\sigma}(0) = 0$ for all $\sigma$ of length~$d$, and the self-input in~(\ref{eqn:ones-reg-lowest}) assures that the same applies to all $t \geq d$, regardless of the trajectory of the variable $i_q$.  Now the lemma follows from our choice of~$d$. $\Box$

\bigskip

Observe that in bi-quadratic Boolean networks, if there exists a variable $x$ whose outdegree is $<2$, then there must exist a variable $y$ whose indegree is $<2$. Therefore, we can add to our system a copy of $B_q$ by defining new regulatory functions that connect $x$ to $i_q$ and $o_q$ to $y$ using the conjunction regulatory functions, and keep repeating this procedure until there are no variables with outdegree $<2$. If there are no such variables left, then, clearly, there are no variables with the indegree $<2$ either.  Moreover, since with probability $\geq q$ any external inputs to any of the copies of~$B_q$ will have no effect on the output of~$B_q$, the arguments in the proof of~\cite{JMI} carry over to the modified system.

It remains to check that we will not add too many variables in this way, and that copies of $B_q$ will start producing value $1$ at $o_q$ sufficiently fast (so they don't affect the workings of  $\bB$). The total number of variables in the Boolean input-output systems $B_1, B_2, B_3$ can be made to scale like $O((\log(n))^{1.5})$.  Unfortunately, this was not explicitly stated in this form in~\cite{JMI}, since for the construction to work, we only needed that the total number of variables in the set $Y = [N] \backslash X$ is bounded from above by a fixed constant times~$(\log(n))^2$ (see~(10) of~\cite{JMI}).  But for $B_3$ the stronger scaling law follows from the formulation of Lemma~10.4 of~\cite{JMI}; for $B_1$ and~$B_2$ it follows from the proofs of Lemmas~5.1 and~5.3 that are given in~\cite{JMI}.

For a fixed $q'=1-x<1$, we need to choose $q<1$ such that with probability $q'$ \emph{each} of the copies of $B_q$ will start generating  the value $1$ at  their variables $o_q$. This will be true for $q \geq (q')^{1/(c\log(n))^{1.5}}$, where $c>0$ is a constant such that $|B_1| + |B_2| + |B_3| \leq (c\log(n))^{1.5}$.
If we let $q=1-\frac{1}{(c\log(n))^{1.5}}x$, then
\[ q^{(c\log(n))^{1.5}} \geq 1-x. \]
For this choice of~$q$, by Lemma \ref{lem:ones}, the number of variables in each copy of $B_q$ is at most
\[ -\log(1-q)=(1.5)c\log(\log(n))-\log(x) \leq c'\log(\log(n))  \]
for some constant $c'>0$, and the total number of new variables will scale like
\[ O((\log(n))^{1.5}\log(\log(n))), \]
which is in compliance with~(10) of~\cite{JMI}.

Now the modified system will, with probability $> q'$, work exactly like the original system for all times $t > d$, where $d = O(\log(\log\log((n)))$ is the depth of the circuit~$B_q$.  The first $d$ steps where $B_1$ and $B_2$ may work improperly have only the effect of slightly increasing the size of the set~$Q$, but not by an order of magnitude (recall that~$Q$ was the set of those indices~$i$ for which the corresponding~$X_i$ was the output of~$B_1$ before everything started working properly). 
 
 With the input-output system~$B_3$ we need to proceed somewhat more carefully, since it needs to work properly right from the beginning.  Recall that in the construction of $B_3$ given in Section~10.3 of~\cite{JMI}, the regulatory functions at the lowest level were already strictly quadratic.  The next levels were designed to produce an ordered version $\vs_o$ of the output $\vs_\ell$ of the lowest level, with $|\vs_o| = |\vs_\ell|$ and all zeros in $\vs_o$ preceding all ones. We quoted a construction from~\cite{AKS} and it is not clear from the quoted result whether this part of the system is strictly quadratic. 
  
The number of such variables in $B_3$ that might take only a single input does not exceed the overall size of~$B_3$, which is $O((\log(n))^\alpha)$ for every $\alpha > 1$  by
the paragraph preceding Lemma~10.4 of~\cite{JMI}.  Thus we can add $O(\log(n)^{1.1})$ copies of the circuits $B_q$ of depth $d_q = O(\log(N)) = O(\log(\log(n)))$ with a total of $O((\log(n))^{1.5})$ variables to give second inputs to these variables.  Moreover, we can  add~$d_q$ levels above the lowest one and use the strictly bi-quadratic Boolean circuit~$B_{c,r}$ for robust coding to produce a version of the lowest level that will preserve its size and will be available for further processing by the original system once all the variables~$o_q$  that are to be used in modifying the subsequent levels have reached their target value~1.
This again requires adding at most $O(|B_3| \log(N)) =O((\log(n))^{1.5})$ new variables and does not violate our restrictions on the size of the set~$Y$ of variables outside the union of the sets $X_i$.

In order to get systems of size exactly~$N$ for all sufficiently large~$N$ we may need to add also some dummy variables (see~\cite{JMI} for an estimate of the size of this set), but these can simply be connected among themselves with AND gates and they don't have any influence on the overall dynamics.

In particular, using subtle coding together with the modifications outlined in the last few paragraphs gives us networks that satisfy the conclusion of Theorem~\ref{thm:tenroot}.  $\Box$

\bigskip

\noindent
\textbf{Proof of Theorems~\ref{thm:cppth} and~\ref{thm:sqrt3}:} Fix $0 < p < 1 < c < 2$, and an auxiliary constant $c_1$ with $c<c_1<2$. For the proof of Theorem~\ref{thm:sqrt3} we make the more stringent assumption that $c < c_1 < \sqrt{3}$. We will show that as long as $N$ is sufficiently large, there exists a  $p$-$c$-chaotic and $p$-unstable $N$-dimensional bi-quadratic cooperative Boolean network. Our strategy will be to first choose some $N_1 < N$ and an $N_1$-dimensional Boolean system $\bB_1 = (2^{N_1}, f)$ that is $p$-$c_1$-chaotic.  We can assume that $\bB_1$ has been constructed as above and has properties (P1)--(P5). For the proof of Theorem~\ref{thm:sqrt3} we will assume in addition that $\bB_1$ is strictly bi-quadratic and uses the robust coding scheme.  We will construct an extension $\bB = (2^N, g)$ of $\bB_1$ so that $g_j = f_j$ for all $j \in [N_1]$ and no variable in $[N_1]$ takes input from any variable in $[N] \backslash [N_1]$. In the proof of Theorem~\ref{thm:sqrt3} we will make an exception for variables~$i_q$  of some copies of $B_q$ of Lemma~\ref{lem:ones} that will receive a second input from $[N] \backslash [N_1]$.  This  provision will preserve the property of $p$-$c_1$-chaos in $\bB_1$, in the sense that with probability~$> p$ a randomly chosen trajectory will reach an attractor of length~$> c_1^{N_1}$.  Note that this implies $p$-$c$-chaos in~$\bB$ as long as $c_1^{N_1} \geq c^N$, or, equivalently,

\begin{equation}\label{eqn:ratio-N/N1-first}
\frac{N_1}{N} \geq \frac{\ln c}{\ln c_1}.
\end{equation}

For a given $N$, let $N_1= \lceil \frac{\ln c}{\ln c_1}N \rceil$. We will first present a construction of the extension $\bB$ of $\bB_1$. At the end of the proof we will argue that for sufficiently large~$N$ the number of required new variables is sufficiently small so that there is enough room for them in~$[N] \backslash [N_1]$.
The regulatory functions for the variables in $[N] \backslash [N_1]$ will be chosen in such a way that the system detects and keeps a permanent record of a proportion of $> p$ of all single-bit flips in the initial conditions.  Having a ratio $\frac{N_1}{N} \approx 1$ again is very helpful here, since it assures that most of these single-bit flips will occur at variables in~$N_1$, and property~(P5) in turn implies that we may restrict our attention to those single-bit flips that happen at some variable $j \in X_i$ for $i > i_2$.  Notice that any such single-bit flip changes $|s_{X_i}(0)|$ for some $i > i_2$.  By Property~(P2), which is preserved under robust coding, and  by induction we will have
$|s_{X_i}(0)| = |s_{X_{i_2+1}}(i - i_2 -1)|$, which allows us to construct the extension in such a way that the only variables in~$N_1$ that send input to any of the variables in $[N] \backslash [N_1]$ are the ones in~$X_{i_2+1}$.  If $\bB_1$ is based on the original construction in~\cite{JMI}, Property~(P3) allows us to copy $s_{X_{i_2+1}}(t)$ to a Boolean vector $s_P(t+1)$ whose set of variables $P$ is contained in $[N] \backslash [N_1]$.  For the proof of Theorem~\ref{thm:sqrt3} we need to assume that~$\bB_1$ is strictly bi-quadratic and has been constructed as in the proof of Theorem~\ref{thm:tenroot}, but with robust instead of subtle coding.  Recall that in this construction each variable in $X_{i_2+1}$ acts as a second input to a variable~$i_q$ at the lowest level of some copy of~$B_q$.  We need to change these outputs to variables in $[N] \backslash [N_1]$ and reassign new  second input variables from $[N] \backslash [N_1]$ to the newly orphaned variables~$i_q$. As we already know, this operation is not expected to alter the relevant  dynamical properties of~$\bB_1$.

Let $u = u(p)$ be a fixed positive integer whose meaning will become apparent shortly. Now we can incorporate a Boolean circuit $B_4$ into $\bB_1$ whose set of variables is contained in $[N] \backslash [N_1]$ that takes $s_P$ as input, produces $u$ copies
of it  and writes its output to vectors $s_O^w$ for $w \in [u]$ of the same dimension after $d_4$ steps so that for all $w \in [u]$ we have $|s_{X_{i_2+1}}(t)| = |s_P(t+1)| = |s_O^w(t+1 + d_4)|$ and all zeros in $|s_O^w(t+1 + d_4)|$ precede all ones in this vector.  We already know from Proposition~10.1 and the proof of Lemma~5.2 of~\cite{JMI} that this can be accomplished by a cooperative bi-quadratic Boolean circuit of depth $d \leq \gamma_4 \log(m)$ that contains a total of $\gamma_4 m \log(m)$ variables, for some constant $\gamma_4$ that is independent of~$N$, where $m$ is on the order of~$\sqrt{N}$.

We would like to create and keep a permanent record of the values of $s_O^w(t+1 + d_4)$ for all times $t < |I| - i_2$.  If this can be done, then the permanent record will persist throughout the attractor, which implies that every single-bit flip in an initial condition that happens at some variable $j \in X_i$ for $i > i_2$ will move the system to a different attractor.  We would like to keep this record in circular data tapes of Boolean vectors $(s_{Z_i^w}: \, i \in I)$ so that for~$t \geq |I| + d_4$ and for all $w \in [u]$ we have
$s_{Z_i^w}(t+1) = s_{Z_{i+1}^w}(t)$ and also $s_{Z_{|I|-1}^w}(t+1) = s_{Z_{0}^w}(t)$, with the tape holding copies of $s_O^w(t+1+d_4)$ in
$s_{Z_{i+t}^w}(|I|+ d_4)$ at time $|I|+d_4$ for all $1 \leq t \leq |I|$.

There are several technical problems with implementing this idea in its original form.
First of all, in order to not use too many new variables, we will actually record only a part of the values of variables $o_1, \ldots, o_m$ from $O$. For the time being, let us just say that we will choose some $j,J$ with $1<j<m/2<J<m$, and keep track of variables $o_j, \ldots, o_J$ only. The values of $j, J$ will be selected in such a way that the difference between $s(0)$ and its single-bit flip $s^*(0)$ will be visible with sufficiently high probability in the window $o_j, o_{j+1},..., o_J$ after placing all zeros before all ones in $s(0), s^*(0)$.

In order to record anything in a circular data tape, for some $i^*$ the variables in $Z_{i^*}^w$ need to take a second input from variables in $O$ in addition to the input from $Z_{i^*+1}^w$ that will be responsible for the copying of the tape.  For our accounting to work as specified above, we need $i^* = i_2 + 1$, but it will be more convenient to write $i^*$.

This leads to our first technical problem: we need to make sure that the relevant data that has been transferred (regardless of how it is done) from $O$ to $Z_{i^*}^w$ at times $d_4 < t \leq d_4 + |I|$ are not erased at subsequent times.
Let us for the time being assume for simplicity that $u=1$, which allows us to drop confusing superscripts~$w$; the solution to the first problem has a straightforward generalization to $u \geq 1$. Enumerate the variables in~$Z_{i^*}$ as $z_j, z_{j+1}, \ldots, z_J$, the variables in~$Z_{i^*+1}$ as $z_j^+, z_{j+1}^+, \ldots, z_J^+$ and define:

\begin{equation}\label{eqn:Zi*-update}
\begin{split}
s_{z_\mu}(t+1) &= s_{z^+_\mu}(t) \vee s_{o_\mu}(t) \quad \mbox{for} \quad j \leq \mu \leq m/2,\\
s_{z_\mu}(t+1) &= s_{z^+_\mu}(t) \wedge s_{o_\mu}(t) \quad \mbox{for} \quad m/2 < \mu \leq J.
\end{split}
\end{equation}

This definition assures that if $|s_P| = m/2$, which is true for all times $t \geq |I|$, then the input from the variables in~$O$ has no effect whatsoever, since in this case the first half of the variables of~$O$ evaluate to~0 and the second half to~1. In particular, by property~(P4a) this will be the case, with probability~$> p$, whenever~$O$ records the size of a vector $s_{X_i}(k|I|)$ for some $k > 1$ with $i > i_2$.

Unfortunately, this definition does not warranty that \emph{exact} copies of $s_O$ will be transferred to $Z_{i*}$. 
Let us focus on the case where $j \leq \mu \leq m/2$; the argument for the case of $\mu > m/2$ is dual.
Each value $s_{o_\mu}(t) =1$ gets faithfully copied to $s_{z_\mu}(t+1)=1$, but the updating as specified by~(\ref{eqn:Zi*-update}) will also introduce some random occurrences of $s_{z_\mu}(t+1)=1$ while $s_{o_\mu}(t)=0$, due to $s_{z^+_\mu}(t)=1$. But consider a situation where the $s_{z^+_\mu}(t)$ are random and we want to use~(\ref{eqn:Zi*-update}) to record to the data tape the sizes of some~$s_{X_i}(0), s^*_{X_i}(0)$ for $i > i_2$, as coded by the Boolean variables~$s_{o_\mu}(t)$, that differ by a single-bit flip  and are  such that $j \leq |s_{X_i}(0)| = \mu < \mu+1 = |s^*_{X_i}(0)| \leq J$.   
 
 However, as long as $s_{z^+_\mu}(t)=0$, a~1 will be copied to~$s_{z_\mu}(t+1)$ only for the trajectory of the corresponding $\vst(0)$, but not for the trajectory of
$\vs(0)$.  If this happens, the two trajectories will reach different attractors and we will say that our recording tape \emph{successfully distinguishes these two initial conditions.} In the proof of Theorem~\ref{thm:cppth} we can assume that $s_{z^+_\mu}(t)$ takes the value~$0$ with probability~$0.5$, which therefore is the probability that a given recording tape will successfully distinguish the two initial conditions as specified.  These events are independent for the~$u$ data tapes, thus by choosing $u$ large enough so that $0.5^u < 1 -p$ we can assure that the probability of success in at least one recording tape is~$> p$, which is all we need for $p$-instability.

The third problem we need to take care of is to choose the values of $j,J$.
We need that with probability~$>p$ we will have

\begin{equation}\label{eqn:CTL-sizeSXi}
j \leq |s_O^w(t+1 + d_4)| \leq J
\end{equation}
as long as $s_O^w(t+1 + d_4)$ records the size of some $s_{X_i}(0)$,
so that a random single-bit flip in this vector  can alter the permanent record.
 By the Central Limit Theorem and Chebysheff's Inequality, this can be achieved, for sufficiently large~$N$ and hence~$m$, if
 $j \leq m/2 - \gamma_5 \sqrt{m}$ and $J \geq m/2 + \gamma_5 \sqrt{m}$ for some constant~$\gamma_5$ that depends on~$p$, but not on~$m$.
Since~$u$ does not depend on $N$ and $m = O(\sqrt{N})$, using such $j, J$ we will be able to construct recording tapes that altogether use only on the order of $u\sqrt{m}|I|$ or $N_1^{3/4}$ variables. For the proof of Theorem~\ref{thm:sqrt3} we will also need $O(N_1^{3/4})$ copies of the Boolean input-output system $B_q$, where $q$ can be chosen as  $1 - \frac{1}{N_1}$. This will add another $O(N_1^{3/4}\log(N_1))$ variables. By connecting these as in the proof of Theorem~\ref{thm:tenroot} to achieve a strictly bi-quadratic network. The modification of all monic regulatory functions of the circuit $B_4$ and the data record tapes may result in missing a few single-bit flips in $s_{X_i}(0)$ for $i_2 < i < i_2 + d$ that $\bB_1$ detects, where $d$ is the depth of $B_q$, but this is no problem, since $d = O(\log(N_1))$ is very small relative to~$|I|$. In either case, the total number  $M$ of variables  in $\bB$ does not exceed $N$, provided that $N$ is sufficiently large. If $M<N$, we add to $\bB$ some dummy variables and connect them as in the proof of Theorem~\ref{thm:tenroot}. These dummy variables will not destroy $p$-instability of $\bB$.
$\Box$

\bigskip

\noindent
\textbf{Proof of Theorems~\ref{thm:decoh} and~\ref{thm:decoh-bi}:} For the proof of Theorem~\ref{thm:decoh}, fix $0 < \alpha < 0.5$, $0 < p < 1 - 2\alpha < 1 < c < 2^{1-2\alpha}$; for the proof of Theorem~\ref{thm:decoh-bi}, fix  $0 < p < 1  < c < \sqrt{3}$.  We will show that as long as $N$ is sufficiently large,
 there exists a  $p$-$c$-chaotic  $N$-dimensional (strictly) bi-quadratic cooperative Boolean network with the required decoherence property. Similarly to the previous proof, our strategy will be to first choose some $N_1 < N$  and an $N_1$-dimensional Boolean system $\bB_1 = (2^{N_1}, f)$ that is constructed as in the previous proof and satisfies the conclusion of Theorem~\ref{thm:cppth} (in the case of the proof of Theorem~\ref{thm:decoh}) or Theorem~\ref{thm:sqrt3} (in the case of the proof of Theorem~\ref{thm:decoh-bi}) for some auxiliary constants~$p_1$ and~$c_1$.
 For Theorem~\ref{thm:decoh} we will choose $p_1 = \frac{p}{1-2\alpha}$ and $c_1 = c^{1/(1-2\alpha)}$ and
 for Theorem~\ref{thm:decoh-bi} we will choose any  $p_1$ with $p < p_1 < 1$ and $c_1$ with $c < c_1 < \sqrt{3}$.

In the proof of Theorem~\ref{thm:decoh} we will construct an extension $\bB = (2^{N}, g)$ of $\bB$ so that $g_j = f_j$ for all $j \in [N_1]$ and no variable in $[N_1]$ takes input from any variable in $[N] \backslash [N_1]$. This latter provision will preserve the property of $p_1$-$c_1$-chaos in $\bB_1$, in the sense that with probability~$> p_1$ a randomly chosen trajectory will reach an attractor of length~$> c^{N_1}$. Moreover, $p_1$-instability of $\bB_1$ will be preserved in the sense that a proportion of $>p_1$ of single-bit flips in initial conditions that occur at variables in $[N_1]$ will result in trajectories that reach different attractors. However, in contrast to the proof of Theorem~\ref{thm:cppth}, we will no longer aim for making the ratio $\frac{N_1}{N}$ arbitrarily close to~1; instead, we will choose $N_1 = \lceil (1-2\alpha)N\rceil$, which gives

\begin{equation}\label{eqn:ratio-N1/N2}
\frac{N_1}{N} \approx 1-2\alpha.
\end{equation}

This has two important consequences that are reflected in the statement of Theorem~\ref{thm:decoh}: First of all, $c_1$-chaos in $\bB_1$ will guarantee at most $c_1^{1-2\alpha}$-chaos in $\bB$, which is the same as $c$-chaos by our choice of~$c_1$.  Second, the proportion of single-bit flips in initial conditions that occur in variables in~$N_1$ is at most $1-2\alpha$. Thus
$p_1$-instability in $\bB_1$ translates at most into $(1-2\alpha)p_1$-instability in~$\bB$, that is, $p$-instability.

In the proof of Theorem~\ref{thm:decoh-bi} we will need to alter some regulatory functions in $\bB_1$ for the variables~$i_q$ of some copies of $B_q$, as we did in the proof of Theorem~\ref{thm:sqrt3}.  By the same argument as in the previous proof,  this is not expected to alter the essential features of the dynamics of $\bB_1$. We will choose $N_1= \lceil \max \left\{\frac{\ln c}{\ln c_1}N, \frac{p}{p_1}N\right\} \rceil$, which ensures, by the same argument as in the proofs of
Theorems~\ref{thm:cppth} and~\ref{thm:sqrt3}, that $\bB$ will be $p$-$c$-chaotic and $p$-unstable.

Thus in both constructions, a random single-bit flip in a randomly chosen initial condition will leave a permanent record in at least one of the data record tapes of~$\bB_1$.
Let $Z_i^w$ be as in the proof of Theorems~\ref{thm:cppth} and~\ref{thm:sqrt3}. Let $z_{\mu}^w$ denote the $\mu$-th element of~$Z_{i_3}^w$.
By the construction in the proof of Theorems~\ref{thm:cppth} and~\ref{thm:sqrt3} and our choice of the ratios $N_1/N, p/p_1$, we can conclude that if $\vs(0), \vst(0)$ are two randomly chosen initial conditions of $\bB$, then with probability~$>p$ the following will hold:

\begin{equation}\label{eqn:detec-mu}
\exists t_0, \mu \in [J-j+1], w \in [u] \forall  k > 0 \ s_{z_{\mu}^w}(t_0 + k|I|) \neq s_{z_{\mu}^w}^*(t_0 + k|I|),
\end{equation}
where $s_i(t), s_i^*(t)$ denote the values of variable~$i$ at time~$t$ in the trajectories of $\vs(0), \vst(0)$ respectively.

The next step in the construction is to add a Boolean input-output system~$B_6$  to $\bB_1$ that copies the values of the variables $z_{\mu}^w$ at selected times to a single variable~$k^*$ so that~(\ref{eqn:detec-mu}) will imply

\begin{equation}\label{eqn:detect-k*}
\forall t \exists t^+ > t \ s_{k^*}(t^+) \neq s_{k^*}^*(t^+).
\end{equation}

Let $(J-j+1)u < T < |I|$ be a prime number.  Since $(J-j+1)u = O(\sqrt{|I|})$,  by the Prime Number Theorem, such $T$ exists for sufficiently large~$N$.
Let $\nu:  ([J] \backslash [j-1]) \times [u] \rightarrow [(J-j+1)u]$ be a bijection. For each  $\mu \in [J-j+1]$ and $w \in [u]$  define a vector $\vec{r}_{\mu, w} \in 2^R$ of length~$|R|=(J-j+1)u$ that takes the value~1 only on its $\nu(\mu,w)$-th coordinate $r_{\nu(\mu,w)}$ and takes the value~0 otherwise.   Lemma~10.4 of~\cite{JMI} implies the existence of a Boolean input-output system $B_{5}$ with output vector~$\vec{r} \in 2^R$ such that with probability arbitrarily close to~1,

\begin{equation}\label{eqn:B_5-workings}
\forall \, (\mu, w)\in ([J]\backslash [j-1]) \times [u] \exists \, t_{\mu, w} \forall \, k \in \NN \ \vec{r}(t_{\mu, w} + kT)  = \vec{r}_{\mu, w}.
\end{equation}

Moreover, $B_{5}$ requires adding only $O((J-j+1)uT\log((J-j+1)uT)) = O(N^{3/4}\log(N))$ variables.

 Create a new set of variables  $R^*$ with $|R^*| = |R|$ and define regulatory functions for the Boolean vector $r^*$ with this domain by

\begin{equation}\label{eqn:reg-r*}
r^*_{\nu(\mu,w)}(t+1) = r_{\nu(\mu,w)} \wedge s_{z_{\mu}^w}(t).
\end{equation}

Make $k^*$ the output of a Boolean circuit~$B_6$ that calculates the conjunction of all the variables in $R^*$.  Since $T$ is relatively prime with $|I|(J-j+1)u)$, this guarantees that the value of each variable in the union of all data tapes will be copied infinitely often to~$k^*$ and gives the implication (\ref{eqn:detec-mu}) $\Rightarrow$ (\ref{eqn:detect-k*}).

Finally, for the proof of Theorem~\ref{thm:decoh} we add another  Boolean circuit $B_7$ to $\bB_1$ that is composed of variables in $[N] \backslash [N_1]$ and copies the value of its single input variable $k^*$ to $\lceil \alpha N \rceil$ distinct output variables after $d_7$ time steps.  By Proposition~10.1 of~\cite{JMI}, there exists a cooperative bi-quadratic Boolean circuit that accomplishes this task and uses at most $2\lceil \alpha N \rceil$ variables. Thus the addition of~$B_7$  does not allow us to achieve a higher ratio $\frac{N_1}{N}$ than in in~(\ref{eqn:ratio-N1/N2}). But since~$B_7$ is the most expensive part of the construction  in terms of the number of necessary additional variables; the ratio can be arbitrarily
close to $1 - 2\alpha$.

Let us recapitulate how this  construction ensures $p$-$\alpha$-decoherence.  A random single-bit flip $\vst(0)$ in a randomly chosen initial condition~$\vs(0)$ happens with probability
$\approx 1- 2\alpha$ at a variable in~$[N_1]$.  By the proof of Theorem~\ref{thm:cppth}, with probability that can be chosen arbitrarily close to~1, it will leave a permanent record in at least one of the data record tapes in~$N_1$. This record will result in infinitely many times~$t^+$ where the trajectories differ at variable~$k^*$, as in~(\ref{eqn:detect-k*}). This difference in turn will be amplified by~$B_7$ to~$\lceil \alpha N \rceil$ distinct variables, and $p$-$\alpha$-decoherence follows. 

For the proof of Theorem~\ref{thm:decoh-bi}, $B_7$ will copy  $k^*$ only  to $\frac{N}{\Theta \log(N)}$ distinct output variables, where $\Theta$ will be determined shortly.  This will ensure $p$-$\frac{N}{\Theta \log(N)}$-decoherence as in the previous paragraph and requires adding at most $\frac{2N}{\Theta \log(N)}$ variables.   The resulting system will not yet be strictly bi-quadratic; we will need to add copies of the Boolean circuits $B_q$ where $1 - q < \frac{1}{N}$  and connect their respective output variables~$o_q$ by conjunctions to variable $j \in B_{6} \cup B_7$ with monic regulatory functions. This requires adding $\frac{2N}{\Theta} + O(N^{3/4}\log(N)\log(\log(N)))$ new variables.  By choosing $\Theta$ large enough so that $\frac{2N}{\Theta} \ll (1-2\alpha)N$ we don't exceed the allotment of $N - N_1$ additional variables specified by~(\ref{eqn:ratio-N1/N2}).  Finally, we  add dummy variables if needed, and  connect variables that have fewer than 2 outputs to variables~$i_q$ of copies of~$B_q$.  The resulting system will have the properties specified in Theorem~\ref{thm:decoh-bi}. $\Box$

\section{Conclusion and future directions}\label{concludesec}

In this paper and its prequel~\cite{JMI} we studied the problem whether cooperativity, that is, the total absence of negative interactions, precludes certain types of chaotic dynamics in Boolean networks, at least under additional assumptions on the number of inputs and outputs per node. This is a natural question in view of the analogous result for continuous flows that was mentioned in the introduction, and the well-known fact that Boolean networks with few inputs per node tend to have ordered dynamics.

Chaotic dynamics of Boolean networks is characterized by very long attractors, very few eventually frozen nodes, and high sensitivity to perturbations of initial conditions.  While these three hallmarks usually go together, the answer to our question crucially depends on how chaos in Boolean networks is formalized.  The notion of $p$-$c$-chaos formalizes genericity of very long attractors and also implies genericity of very few eventually frozen nodes.  We showed that cooperativity does not impose any nontrivial bounds on this property, even in bi-quadratic Boolean networks. Similarly, in strictly bi-quadratic networks, cooperativity does not imply additional bounds on $p$-$c$-chaos beyond the previously known bound of $c < 10^{1/4}$ for $c$-chaos.

However, the situation changes when one considers notions of high sensitivity to perturbations of initial conditions.  The strongest such notion considered here,
$p$-$\alpha$-$q$-decoherence, while still possible in $p$-$c$-chaotic Boolean networks in general, is outright precluded by cooperativity. The weakest of these three notions, $p$-instability, is still consistent with $p$-$c$-chaos in cooperative bi-quadratic Boolean networks for all $0 < p < 1 < c < 2$. But if in addition it is assumed that the network is strictly bi-quadratic, a stronger bound $c < \sqrt{3}$ applies, and the bound is strict.  

The notion of $p$-$D$-decoherence comes in many flavors, depending on the parameter~$D$. While it is consistent with cooperativity and $p$-$c$-chaos for all $0 < p < 1 < c < 2$ and all meaningful linear~$D$ in general, we were only able to construct bi-quadratic and strictly bi-quadratic Boolean systems that satisfy this property under some additional restrictions on~$p$, $c$, and~$D$.  Since any form of $p$-$D$-decoherence implies $p$-instability, there must be some such restrictions at least under the additional assumption that the system is strictly bi-quadratic. We conjecture that there are some restrictions for bi-quadratic systems as well.  However, it remains an open problem to find the optimal upper bounds on the amount of $p$-$c$-chaos and $p$-$D$-decoherence that can simultaneously occur in such networks.

Thus cooperativity, by itself and in conjunction with suitable restrictions on the interactions of the variables, does impose restrictions on how much chaos is possible in a Boolean network. These results can be interpreted as counterparts of the corresponding theorem for flows.  Our work shows that valid results of this kind require very specific conceptualization of certain hallmarks of chaotic dynamics.  We believe that these subtleties  need to be well understood if researchers are to make valid inferences from dynamical properties of a Boolean approximation to an ODE model of a natural system about the ODE dynamics or the behavior of the natural system itself.

\end{document}